\documentclass[12pt,twoside,a4paper]{amsart}
\usepackage{t1enc,amssymb}
\usepackage{times,a4wide}

\newtheorem{theorem}{Theorem}
\newtheorem*{thm2prime}{Theorem 2'}

\newtheorem{lemma}[theorem]{Lemma}

\newcommand{\la}{\langle}
\newcommand{\ra}{\rangle}

\newcommand{\Z}{\mathbb{Z}}

\newcommand{\R}{\mathbb{R}}
\newcommand{\T}{\mathbb{T}}

\begin{document}
\title{$L^\infty$ to $L^p$ constants for Riesz projections}
\author{Jordi Marzo}
\author{Kristian Seip}
\address{Department of Mathematical Sciences,
Norwegian University of Science and Technology, N-7491 Trondheim, Norway}
\email{jordi.marzo@math.ntnu.no}
\address{Department of Mathematical Sciences,
Norwegian University of Science and Technology, N-7491 Trondheim,
Norway} \email{kristian.seip@math.ntnu.no}


\thanks{The authors are supported by the Research Council of
Norway grant 185359/V30.} 

\date{\today}
\begin{abstract}
The norm of the Riesz projection from $L^\infty(\T^n)$ to
$L^p(\T^n)$ is considered. It is shown that for $n=1$, the norm
equals $1$ if and only if $p\le 4$ and that the norm behaves
asymptotically as $p/(\pi e)$ when $p\to \infty$. The critical
exponent $p_n$ is the supremum of those $p$ for which the norm
equals $1$. It is proved that $2+2/(2^n-1)\le p_n <4$ for $n>1$; it
is unknown whether the critical exponent for $n=\infty$ exceeds $2$.
\end{abstract}

\maketitle

\section{Introduction}

This work originated in an attempt to answer the following question:
Do there exist pairs of exponents $q$ and $p$ with $q>p>2$ for which
the Riesz projection on the infinite torus is bounded from $L^q$ to
$L^p$? This question remains open, as far as we know. The present
note presents a few results of some intrinsic interest in the finite
dimensional setting, giving relevant background for the original
problem for the infinite torus.

Using standard multi-index notation, we write the Fourier series of
a function $f$ in $L^2(\T^n)$ on the $n$-torus $\T^n$ as
\[ f(\zeta) = \sum_{\alpha \in \Z^n} \hat{f}(\alpha)\zeta^{\alpha}.
\]
The operator
\[ P^+_n f(\zeta)= \sum_{\alpha \in \Z_+^n}
\hat{f}(\alpha)\zeta^{\alpha} \] is the Riesz projection on $\T^n$,
and, as an operator on $L^2(\T^n)$, it has norm $1$. If we instead
view $P_n^+$ as an operator on $L^p(\T^n)$ for $1<p<\infty$, then a
theorem of B. Hollenbeck and I.~Verbitsky \cite{HV00} says that its
norm is $(\sin \frac{\pi}{p})^{-n}$.

We compute the norm $\|f\|_p$ of a function $f$ in $L^p(\T^n)$ with
respect to Lebesgue measure $\sigma_n$ on $\T^n$, normalized such
that $\sigma_n(\T^n)=1$. Using this normalization, we let
$\|P_n^+\|_{q,p}$ denote the norm of the operator $P_n^+: L^q(\T^n)
\to L^p(\T^n)$ for $q\ge p\ge 2$. We will restrict ourselves to
computations and estimates of the norms $\|P_n^+\|_{\infty,p}$. By
H\"{o}lder's inequality, $p\mapsto \|P_n^+\|_{\infty,p}$ is a
continuous and nondecreasing function, and, by the theorem of
Hollenbeck and Verbitsky, we have $\|P_n^+\|_{\infty, p} \le (\sin
\frac{\pi}{p})^{-n}$. Of particular interest is the number
\[ p_n = \sup\left\{p\ge 2:\ \|P_n^+\|_{\infty,p}=1\right\}, \]
called the critical exponent of $P_n^+$ according to the terminology
of \cite{FIP84}. The critical exponent $p_n$ is well-defined since
clearly $\|P_n^+\|_{\infty,2}=1$. By continuity, we have
$\|P_n^+\|_{\infty,p_n}=1$.

We will present three theorems. The first says that the critical
exponent of $P_1^+$ equals $4$. In view of this result, one is led
to ask if the precise value of $\|P_1^+\|_{\infty,p}$ can be
computed also when $p>4$ and whether we can compute or estimate the
critical exponent $p_n$ for $n>1$. These problems are only given
partial solutions: Our second theorem gives the right asymptotics
for $\|P_1^+\|_{\infty,p}$ when $p \to \infty$, and our third
theorem says that $2+2/(2^n-1)\le p_n < 4$ when $n>1$.

The next three sections present these results. Section 5 contains a
brief discussion of the problem for the infinite torus, while the
final section discusses extensions to the setting of compact abelian
groups.

\section{The critical exponent of $P_1^+$}

\begin{theorem}                                        \label{theo1}
The critical exponent of $P_1^+$ is $4$.
\end{theorem}

\proof We write $P_1^-=I-P_1^+$ and note that $(P_1^+f)^2\perp
(P_1^-f)^2$ whenever $f$ is a bounded function on $\T$. Thus
\[ \| P_1^+ f \|_4^4  =\| (P_1^+f)^2 \|_2^2\le \| (P_1^+f)^2-(P_1^-f)^2 \|_2^2=\| f(P_1^+ f -P_1^-f)
\|_2^2
    \le \| f \|_\infty^2 \| f \|_2^2. \]
This estimate implies that $p_1\ge 4$. To see that we also have
$p_1\le 4$, we consider the function $f(\zeta)=(1-\epsilon
\zeta)^2/|1-\epsilon \zeta|^2$. We assume that $0<\epsilon<1/2$ and
find that $P_1^+ f(\zeta)= 1-\epsilon^2 - \epsilon \zeta.$ We
estimate the $L^p$ norm of $P_1^+f$ from the power series expansion
of $(1-\epsilon \zeta/(1-\epsilon^2))^{p/2}$. This leads to the
estimate
\[ \| P_1^+f \|_p^p=1+\left(\frac{p^2}{4}-p\right) \epsilon^2 +
O(\epsilon^4) \] when $\epsilon \to 0$. It follows that we may
achieve $\| P_1^+ f\|_p>1$ for every $p>4$ by choosing $\epsilon$
sufficiently small.\qed

We note that in general we have $(P_n^+f)^2\perp ((I-P_n^+)f)^2$
only when $n=1$, so that the preceding proof does not work when
$n>1$.

\section{Asymptotic behavior of $\|P_1^+\|_{\infty,p}$ when $p\to\infty$}

\begin{theorem}
   We have $\lim_{p\to \infty}p^{-1} \| P_1^+ \|_{\infty,p}=(\pi e)^{-1}.$
\end{theorem}

This theorem is a corollary of a corresponding result for the
Hilbert transform (the conjugation operator), which we define as
    $$Hf(\zeta)=\tilde{f}(\zeta)=-i \sum_{k\in \Z} \mbox{sign}(k) \hat{f}(k)\zeta^{k}.$$
By a well-known theorem of Pichorides \cite{Pic72}, we have
    $$ \| H \|_{p,p}=\max\left\{ \tan \frac{\pi}{2 p},\cot \frac{\pi}{2 p}\right\}.$$
The Hilbert transform maps real functions to real functions, and we
write $H_{\R}$ when the domain is a real $L^p$ space.

\begin{thm2prime}
We have $\lim_{p\to \infty}p^{-1}\| H \|_{\infty,p}=\lim_{p\to
\infty}p^{-1}\| H_\R \|_{\infty,p}=2(\pi e)^{-1}.$
\end{thm2prime}
    The following result of Zygmund \cite[Theorem 2.11, chap. VII, vol. 1]{Zyg68} will give an upper bound
    for $\| H_\R \|_{\infty,p}.$
\begin{lemma}(Zygmund)
    For real valued $f$ such that $|f|\le 1$ and $0\le
    \alpha<\pi/2$, we have
    $$\frac{1}{2\pi}\int_0^{2\pi}e^{\alpha |\tilde{f}(e^{i\theta})| }d\theta\le \frac{2}{\cos \alpha}.$$
\end{lemma}

\begin{proof}[Proof of Theorem 2']
    From Zygmund's theorem and Chebychev's inequality, we get that for real valued
    $f$ with $|f|\le 1$ we have
    $$\sigma_1\left( \{ \zeta  :\ |\tilde{f}(\zeta)|> \lambda \}\right)  \le
    \frac{2}{\cos \alpha}e^{-\alpha \lambda},$$
and thus
\[
    \|\tilde{f}\|_p^p =p\int_0^{\infty} \lambda^{p-1} \sigma_1\left(
    \{ \zeta  : \ |\tilde{f}(\zeta)|> \lambda \}\right)  d\lambda
\le \frac{2p }{\alpha^p \cos \alpha} \Gamma(p+1). \] Now Stirling's
formula implies that
$$\lim_{p\to \infty} \frac{1}{p}\| H_\R \|_{\infty,p}\le  \frac{2}{\pi e}.$$

To prove the reverse inequality,  we consider the function
    $$ f(e^{i\theta})=\arg(1-e^{i\theta})= \left\{\begin{array}{ccc} \theta/2-\pi/2, &  & 0\le \theta\le \pi \\
                                                    \theta/2+\pi/2 , & &
                                                    -\pi<\theta<0.
                                                    \end{array} \right.$$
We have $H(\arg(1-e^{i\theta}))=-\log |1-e^{i\theta}|$ and
    $$\| Hf \|_p\le \| H_\R \|_{\infty,p} \| \arg(1-e^{i\theta}) \|_\infty=\| H_\R \|_{\infty,p} \frac{\pi}{2}.$$
Since
\[
    \| \log |1-e^{i\theta}| \|_p^p  =
    \frac{1}{\pi} \int_0^{\pi} \left|\log( 2\sin\frac{\theta}{2} )\right|^p
    d\theta \ge \frac{1}{\pi} \int_0^{1} |\log \theta |^p d\theta
    =\frac{1}{\pi}\Gamma(p+1),\]
it follows that
    $$\frac{2}{\pi e}=\lim_{p\to \infty}\frac{2}{\pi p}\| Hf  \|_p
    \le \lim_{p\to \infty}\frac{1}{p}\| H_\R \|_{\infty,p},$$
    and we conclude that
    $\lim_{p\to \infty} p^{-1}\| H_\R \|_{\infty,p}=2(e \pi)^{-1}$.

We turn next to the complex case. What follows is a small variation
of a construction used to prove vector-valued inequalities \cite[pp.
311--315]{Gra04}. We begin by noting that
    for arbitrary real numbers $\lambda_1$ and $\lambda_2$ and
    $0<p<\infty$, we have
    $$\int_{\R^2} |x_1 \lambda_1 + x_2 \lambda_2|^p e^{-\pi |x|^2}dx=A_p^{p}(\lambda_1^2+\lambda_2^2)^{p/2},$$
    where
    $$A_p=\left( \frac{\Gamma\left( \frac{p+1}{2}\right)}{\pi^{\frac{p+1}{2}}} \right)^{1/p}.$$
    if $f=f_1+i f_2$ is complex valued function with $f_1$ and $f_2$ are real valued, then
\begin{align*}
    \| Hf \|_p^p & =\frac{1}{2\pi}\int_\T |(Hf_1)^2+(Hf_2)^2|^{p/2}
    d\theta
    \\
    &
    = \frac{A_p^{-p}}{2\pi}\int_\T \int_{\R^2} |x_1 Hf_1 + x_2 Hf_2|^p e^{-\pi |x|^2}dx d\theta
    \\
    &
    = A_p^{-p}\int_{\R^2}  \frac{1}{2\pi}\int_\T  |H( x_1  f_1 + x_2 f_2)|^p d\theta  e^{-\pi |x|^2}dx
    \\
    &
    \le
    A_p^{-p}\int_{\R^2}  (\| H_\R \|_{\infty,p})^p \| x_1  f_1 + x_2 f_2 \|_\infty^p e^{-\pi |x|^2}dx.
\end{align*}
Since $x_1 f_1+x_2 f_2 =$Re$[(x_1-ix_2)(f_1+f_2)]$, we get
$$\|Hf\|_p^p \le A_p^{-p}
(\| H_\R \|_{\infty,p})^p \| f \|_\infty^p  \int_{\R^2} |x|^p
e^{-\pi |x|^2}dx$$
    and therefore
$$\| H \|_{\infty,p} \le  A_p^{-1} \left( \int_{\R^2} |x|^p  e^{-\pi |x|^2}dx
\right)^{1/p} \| H_\R \|_{\infty,p} = A_p^{-1}\pi^{-p/2}
\Gamma(p/2+1) \| H_\R \|_{\infty,p}.$$ Using again Stirling's
formula, we obtain
$$ \lim_{p\to \infty} \frac{1}{p}\| H \|_{\infty,p} \le \lim_{p\to \infty} \frac{1}{p}
\| H_\R \|_{\infty,p}=\frac{2}{\pi e}.$$ Since obviously $\| H_\R
\|_{\infty,p}\le \| H \|_{\infty,p}$, we get the desired result.
\end{proof}
Since
    $$P_1^+ f(e^{i\theta})=\frac{1}{2}(f(e^{i\theta})+i\tilde{f} (e^{i\theta}))+\frac{1}{2}\hat{f}(0),$$
we have
    $$\frac{\| H \|_{\infty,p}}{2}-1\le \| P_1^+ \|_{\infty,p}\le 1+\frac{\| H \|_{\infty,p}}{2},$$
    and we therefore obtain Theorem 2 as an immediate consequence of Theorem 2'.

We note that since $\|P_n^+\|_{\infty,p}\le \|P_n^+\|_{p,p}$,
Hollenbeck and Verbitsky's theorem gives
\[ \limsup_{p\to \infty} p^{-n} \|P_n^+\|_{\infty,p} \le \pi^{-n}. \]
On the other hand, using the function $f(\zeta_1)\cdots f(\zeta_n)$
with $f$ as in the proof of Theorem~2', we obtain
\[ \liminf_{p\to \infty} p^{-n}\|P_n^+\|_{\infty,p} \ge (\pi e)^{-n}. \]
A. Chang and R. Fefferman's counterpart to the John--Nirenberg
theorem (see \cite{CF80}) could be used in place of Zygmund's lemma.
However, we are not aware of any version of the John--Nirenberg
theorem for $\T^n$ that is sufficiently precise to improve our
asymptotic estimates for $n>1$.

\section{Critical exponents for $n>1$}

\begin{theorem}                                                         \label{riez-thorin}
    We have  $ 2+\frac{2}{2^n-1} \le p_n <4 $ when $n>1$.
\end{theorem}

Here the right inequality is of interest because it shows that
$p_1>p_2$, and this means that the problem is truly
multi-dimensional in contrast to the one for the $L^p$ to $L^p$
constants. The left inequality is probably far from optimal; the
main point of this estimate is the fact that we have $p_n>2$ for
every $n$. It appears to be a difficult problem to decide whether
$\lim_{n\to\infty} p_n>2$.

\begin{proof}[Proof of Theorem~\ref{riez-thorin}]
We prove the left inequality by induction on $n.$ By Theorem~1, this
inequality is in fact an equality when $n=1$.

Suppose now that we have
    $\| P_{n-1}^+ \|_{\infty,q}=1$
for $q=2+2/(2^{n-1}-1)$. Consider $P^+_n$ as the composition
of the projections $P^+_{n-1}$ acting on the first $n-1$ variables
and $P_1^+$ acting on the $n$-th variable. If we write $\zeta=(\xi,
\zeta_n)$ with $\xi=(\zeta_1,...,\zeta_{n-1})$, we may write this as
\[ P^+_n f(\zeta)=P_{1,\zeta_n}^+P^{+}_{n-1,\xi} f(\xi,\zeta_n).\]
We now observe that since $\|P^+_1\|_{2,2}=\|P_1^+\|_{\infty,4}=1$,
the Riesz--Thorin theorem implies that we also have $\|
P^+_1\|_{q,p}=1$ when
$$2<q<\infty\;\;\;\mbox{and} \;\;\;p=\frac{4q}{2+q}.$$
Setting \[ p=\frac{4q}{2+q}= 2+\frac{2}{2^n-1},\] we therefore
obtain
\begin{align*}
    \| P_n^+ f \|_p^p & =
    \int_{\T^{n-1}}\int_{\T}|P^+_{1,\zeta_n}P_{n-1,\xi}^+f(\xi,\zeta_n)|^p d\sigma_1(\zeta_n) d\sigma_{n-1}(\xi)
    \\ & \le \int_{\T^{n-1}}\left(
    \int_{\T}|P_{n-1,\xi}^+f(\xi,\zeta_n)|^q
    d\sigma_1(\zeta_n)\right)^{p/q}d\sigma_{n-1}
    \\
    &
        \le
    \left( \int_{\T}   \int_{\T^{n-1}} |P^+_{n-1,\xi}f(\xi,\zeta_n)|^p d\sigma_{n-1}(\xi)
    d\sigma_{1}(\zeta_n) \right)^{q/p} \\
 &   \le \left( \int_{\T}  \left( \sup_{\xi\in \T^{n-1}} |f(\xi,\zeta_n)|\right)^p d\sigma_1(\zeta_n)   \right)^{q/p}
    \le
    \|f\|_\infty^q.
\end{align*}

We clearly have $p_{n+1}\le p_n$, so the only remaining task is to
show that $p_2<4.$ Let $g$ be a homogeneous holomorphic polynomial
    on the bidisk with $\| g \|_\infty\le 1.$ We set $h=(1-g)/(1-\overline{g})$ and find that
\[
    P_2^+ h =P_2^+ ( (1-g)(1+\bar{g}+\bar{g}^2+\dots )  )
=
     1-P_2^+ ( |g|^2 )-g =  1- \|g\|^2_2  -g.
\]
It follows that
$$( P_2^+ h)^2=(1- \|g\|^2_2)^2-2(1- \|g\|^2_2)g+g^2,$$
    and since the functions $1,g,g^2$ are mutually orthogonal, we
may compute $\|P_2^+h\|_4$ explicitly:
\begin{align*}
    \| P^+_2 h \|^4_4  =\| ( P^+_2 h)^2 \|^2_2=1+\| g \|_2^8+\| g \|_4^4-2\| g
    \|_2^4.
\end{align*}
Thus $\|P_2^+h\|_4>1$ whenever    $\| g \|_4^4 \ge 2 \| g \|_2^4.$
One can take, for example, the polynomial
   \[ g(z_1,z_2)=\frac{(z_1+z_2)^3}{10}.\]
\end{proof}

Note that we may obtain a slightly better upper bound by computing
the $L^p$ norm of $P_2^+ h$ from the expansion of $(P_2^+ h)^{p/2}$
into a power series in $g$. By this approach, using the polynomial
\[ g(z_1,z_2)=\frac{(z_1+z_2)^{10}}{1025}, \]
we have found that in fact $p_2\le 3.67632$.

\section{The critical exponent for $n=\infty$}

Let $\sigma_\infty$ denote Haar measure on $\T^\infty$ normalized so
that $\sigma_{\infty}(\T^\infty)=1$, and let $L^p(\T^\infty)$ be the
corresponding $L^p$ spaces. A multi-index $\alpha$ is now a sequence
\[ \alpha=(\alpha_1, \alpha_2,...),
\] where only finitely many of the integers $\alpha_j$ are nonzero. We write
$\alpha\ge \beta$ if we have $\alpha_j\ge \beta_j$ for every $j$.

The Riesz projection of a function $f$ in $L^2(\T^\infty)$ with
Fourier series
\[ f(\zeta)=\sum_{\alpha} \hat{f}(\alpha) \zeta^\alpha, \]
can now be written as
\[ P_\infty^+f(\zeta)=\sum_{\alpha\ge 0} \hat{f}(\alpha)
\zeta^\alpha.\] We define the critical exponent of $P^+_\infty$ as
\[ p_\infty = \sup\left\{p\ge 2:\ \|P_\infty^+\|_{\infty,p}=1\right\}. \]
Note the following difference from the finite-dimensional case: we
have either $\|P^+_{\infty}\|_{\infty,p}=1$ or
$\|P^+_{\infty}\|_{\infty,p}=\infty$, so that
$\|P^+_{\infty}\|_{\infty,p}=\infty$ for $p>p_\infty$.

We want to show that $p_\infty = \lim_{n\to \infty} p_n$. It is
clear that the limit exists and that $p_\infty \le \lim_{n\to
\infty} p_n$. To show that we have equality, we assume that $2\le
p_\infty < \lim_{n\to \infty} p_n$. Let $\varphi$ be a function of
norm $1$ in $L^\infty({\Bbb T}^\infty)$ such that
$\|P_\infty^+\varphi\|_p>1$ for $p_\infty < p < \lim_{n\to \infty}
p_n$. Let $n$ be a positive integer and set
\[ \varphi_n(\zeta_1,...,\zeta_n) =\int_{{\Bbb T}^\infty}
\varphi(\zeta_1,..., \zeta_n, \xi_{n+1}, \xi_{n+2},...)
d\sigma_\infty (\xi).\] Then $\|\varphi_n\|_\infty \le
\|\varphi\|_\infty$. We observe also that
\[ P_n^+\varphi_n(\zeta)=P_\infty^+\varphi(\zeta_1,...,\zeta_n, 0, 0, ...).\]
It is plain that $\|P_n^+\varphi_n\|_p\le \|P_\infty^+\varphi\|_p$.
On the other hand, since $P_n^+\varphi_n\to P_\infty^+\varphi$ in
$L^2$, there is a subsequence $P_{n_k}^+\varphi_{n_k}$ converging to
$P_+\varphi$ almost everywhere. Thus, by Fatou's lemma,
$\|P_\infty^+\varphi\|\le \lim_{k\to \infty}
\|P_{n_k}^+\varphi_{n_k}\|_p$. Hence
$\lim_{n\to\infty}\|P_n^+\varphi_n\|_p=\|P_\infty^+\varphi\|_p$,
which means that $\|P_{n}^+\varphi_n\|_p>1$ for sufficiently large
$n$. This contradicts the assumption that $p<p_n$.

We conclude that if we could prove that $\lim_{n\to\infty} p_n >2$,
then we would have a positive answer to the question asked in the
first paragraph of this note. \vspace{2mm}

\section{Extensions and comments}

The preceding results about critical exponents extend to the
following more general setting. Let $G$ be a compact abelian group
    and let $E$ be a subset of the dual group $\hat{G}.$ We then define the
    $E$-projection of a function $f$ in $L^2(G)$ as
$$P_E f(\omega)=\sum_{\gamma \in E} \hat{f}(\gamma)\la \gamma , \omega \ra,\;\;\omega \in G.$$
    When $E$ generates an order in the dual group $\hat{G}$
    (as it may for connected groups $G$), the proof of Theorem \ref{theo1}
    still works, so that $P_E$ has critical exponent 4.
    Observe also that a direct analogue of Theorem \ref{riez-thorin} can be obtained in this
    case.

    In an attempt to simplify matters, we have studied the following example which
    appears to be the simplest nontrivial case, at least from a computational point of view.
    Take $G=\Z_3$, and consider Riesz projection to be the operator obtained by restricting
    to the set
    $\{ 0,1 \}\subset \Z_3$ in the Fourier domain. The set $\{0,1\}$ does not generate an order, so we
    can not apply the observations made above. However, we may compute the critical exponent
    using the fact that the problem of maximizing the $p$-norm of the projection has the
    following geometrical interpretation. Indeed, we want to compute the maximum of $m_0^p+m_1^p+m_2^p$
    where $m_0,m_1,m_2$ are the lengths of the medians of a triangle with vertices in the closed
    unit disc. It may be seen that the critical exponent equals the solution
    of the equation $2^p+2=3(3/2)^p$, which means that $p_1= 3.08164...$. It is a curious
    fact that this number is also the critical exponent of a different projection in
    \cite[p. 265]{FIP84}.

    The corresponding multi-variable problem seems to be not much easier than the one for $\T^n$.
    Even for $\Z_3^2$ we have not been able to compute the critical exponent
numerically. All we can say is that the critical exponent
    for $P_{\{0,1\}^2}$ is strictly smaller than $p_1$ and in fact $p_2\le 2.93039...$. This is
    far from the corresponding lower bound 2.28107... obtained from the Riesz--Thorin
    theorem, cf. the proof of Theorem \ref{riez-thorin}.

\vskip .2in

\emph{Acknowledgements.} We are grateful to Joaquim
Ortega-Cerd\`{a}, Eero Saksman, Dragan Vukotic, and Brett Wick for
enlightening discussions on the subject matter of this paper.

\end{document}